\documentclass[3p,times,chicago]{elsarticle}

\usepackage{amssymb}

\usepackage{amsmath,amssymb,amsthm,textcomp}
\usepackage{amsfonts,graphicx}
\usepackage[mathscr]{eucal}
\usepackage{color}
\usepackage{hyperref}
\usepackage[table]{xcolor}
\usepackage{booktabs}

\usepackage[]{natbib}
\usepackage{subfig}
\usepackage{epsfig}
\usepackage{caption}
\theoremstyle{definition}
\usepackage{float}
\usepackage{mathtools}  
\usepackage{diffcoeff}  
\floatstyle{plaintop}
\usepackage[figuresright]{rotating}
\usepackage{natbib}
\usepackage{amssymb,amsmath, amsfonts}
\usepackage{enumerate}
\usepackage{amsthm}
\usepackage[utf8]{inputenc}
\usepackage[english]{babel}
\setcitestyle{authoryear,open={(},close={)}} 

\usepackage{graphicx}
\usepackage{caption}
\newtheorem{theorem}{Theorem}[section]

\newtheorem{remark}{Remark}[section]
\newtheorem{lemma}{Lemma}[section]
\newtheorem{definition}{Definition}[section]

\numberwithin{equation}{section}

\newtheorem{example}{Example}[section]
\newtheoremstyle
{remarkstyle}
{}
{11pt}
{}
{}
{\bfseries}
{:}
{     }
{\thmname{#1} \thmnumber{#2} }

\theoremstyle{remarkstyle}

\makeatletter
\def\ps@pprintTitle{%
	\let\@oddhead\@empty
	\let\@evenhead\@empty
	\let\@oddfoot\@empty
	\let\@evenfoot\@oddfoot
}
\makeatother

\begin{document}

\begin{frontmatter}




\title{Generalized Fractional Negative Binomial Process}

\author[label1]{Ritik Soni}
\author[label1]{Ashok Kumar Pathak \corref{cor1}}
\ead{ashokiitb09@gmail.com}
\address[label1]{Department of Mathematics and Statistics, Central University of Punjab, Bathinda, India}
\cortext[cor1]{Corresponding author}

\begin{abstract} 
	In this paper, we introduce a generalized fractional negative binomial process (GFNBP) by time changing the fractional Poisson process with
	 an independent Mittag-Leffler (ML) L\'{e}vy subordinator. We study its  distributional properties and its connection to PDEs. We examine the long-range dependence (LRD) property of the GFNBP and show that it is not infinitely divisible. The space fractional  and the non-homogeneous variants of the GFNBP are explored. Finally, simulated sample paths for the ML L\'{e}vy subordinator and the GFNBP are also presented.
\end{abstract}

\end{frontmatter}


\section{Introduction}
Counting processes and their connections to fractional differential equations have received considerable attention in recent years with wide applications in diverse disciplines of the applied sciences, namely,  physics, image processing, infectious diseases modeling, hydrology, finance, and probability theory (see \cite{timmermann1999multiscale}, \cite{laskin2003fractional}, \cite{kumar2020fractional}, \cite{guler2022forecasting}). The fractional Poisson process (FPP) and the negative binomial (NB) process are the two most commonly used counting processes studied in literature. \cite{laskin2009some} employed the FPPs in defining a new quantum coherent system and also studied  fractional versions of the Bell polynomials, Bell numbers, and the Stirling numbers of the second kind. An application of the FPP in risk theory is discussed by \cite{biard2014fractional}. Some actuarial and clinical trial applications of  NB process have been addressed in  \cite{grandell1997mixed} and \cite{cook2003conditional}.


In recent years, several fractional versions of the NB process have been developed using subordination techniques through stable and inverse stable subordinators (see \cite{beghin2014fractional}, \cite{beghin2015fractional}, and \cite{samy2018fractional}). These generalizations are time-fractional versions of the  NB process. Apart from these, space-fractional versions of the NB process have also been studied in literature (see \cite{orsingher2012space}, \cite{polito2016generalization}, \cite{beghin2018space}, and \cite{maheshwari2023tempered}).
 For $0<\beta<1$, let $\{S_\beta(t)\}_{t \geq 0}$ be a $\beta$-stable subordinator with the Laplace transform $\mathbb{E}[e^{-uS_\beta(t)} ]= e^{-tu^\beta}.$ The inverse $\beta$-stable subordinator  $\{E_\beta(t)\}_{t \geq 0}$  is defined as
\begin{equation*}
E_\beta(t) = \inf\{r \geq0 : S_\beta(r)> t\}, \; t \geq 0.
\end{equation*}
Let $\mu>0$, $\rho>0$, and $\{\Gamma(t)\}_{t\geq 0}$ be a gamma process, where $\Gamma(t)\sim G(\mu, \rho t)$, which denotes the gamma distribution with scale parameter $\mu^{-1}$ and shape parameter $\rho t$. Then NB process can be considered as gamma subordinated variant of the Poisson process.  Recently, \cite{samy2018fractional} presented a fractional NB process (FNBP) $\{ \mathcal{Q}_\beta(t,\lambda)\}_{t\geq 0}$ using gamma subordination in FPP which is characterized as 
\begin{equation*}
\mathcal{Q}_\beta(t,\lambda)  = N_\beta(\Gamma(t), \lambda),
\end{equation*}
where $\{N_\beta(t, \lambda)\}_{t\geq 0}$ is the FPP defined in \cite{meerschaert2011fractional}. The Mittag-Leffler (ML) L\'{e}vy process is a well-known example  of  geometric stable process with non-decreasing paths. It may be used in place of the gamma subordinator to construct several time-changed stochastic processes that  may exhibit additional properties due to  Mittag-Leffler delay. In this paper, we present a generalized fractional negative binomial process (GFNBP)  by time changing  the FPP with
an independent Mittag-Leffler (ML) L\'{e}vy subordinator. This process exhibits overdispersion and long-range dependence (LRD) properties. It is not infinitely divisible and may be useful in different areas.

The paper is structured as follows: In Section 2, we present some preliminary notations and definitions. In Section 3, we define the GFNBP and discuss its  main characteristics along with LRD property. The underlying fractional PDEs are also obtained for the pmf of the GFNBP.  We study the space-fractional and the non-homogeneous version of the GFNBP in Section 4. Finally, we present simulated sample paths for the ML L\'{e}vy subordinator and the GFNBP in Section 5.
\section{Preliminaries} In this section, some notations and definitions are given which will be used in the subsequent sections. Let $\mathbb{R}$ and $\mathbb{C}$ denote the set of real and complex numbers, respectively. Let $\mathbb{Z}_{+}=\mathbb{N}\cup \{0\}$,  where $\mathbb{N}$ is the set of natural numbers.
\subsection{Special functions} Here, we present some special functions which are essential for development of results in this paper.\\
\noindent (i) Three parameters Mittag-Leffler function $L_{\beta, \gamma}^{\alpha}(z)$ is defined as (see \cite{podlubny1999introduction}, \cite{prabhakar1971singular})
\begin{equation}\label{ml123}
L_{\beta, \gamma}^{\alpha}(z) = \sum_{k=0}^{\infty} \frac{z^k}{k!\Gamma(\gamma+\beta k)}\frac{\Gamma(\alpha +k)}{\Gamma(\alpha)}, \;\;\beta, \gamma, \alpha,  z \in \mathbb{C} \text{ and } \text{Re}(\beta)>0, \text{Re}(\gamma)>0, \text{Re}(\alpha)>0.
\end{equation}
\noindent (ii) For $z \in \mathbb{C}$ and $ 0< \alpha <1,$ the M-Wright function $M_\alpha(z)$ is defined by (see \cite{gorenflo2015fractional})
\begin{equation*}
M_\alpha(z) = \sum_{k=0}^{\infty} \frac{(-z)^k}{k! \Gamma (- \alpha n + (1- \alpha))}.
\end{equation*}
The generalized Wright function is defined by (\cite{kilbas2002generalized})
\begin{equation}\label{gwf11}
_p\psi_q \left[z\; \vline \;\begin{matrix}
\left(\alpha_i, \beta_i\right)_{1,p}\\
(a_j,b_j)_{1,q}
\end{matrix} \right] = \sum_{k=0}^{\infty} \frac{z^k}{k!} \frac{\prod_{i=1}^{p} \Gamma(\alpha_i + \beta_i k)}{\prod_{j=1}^{q}\Gamma(a_j + b_j k)},\;\;  z, \alpha_i, a_i \in \mathbb{C}\; \text{and}\; \beta_i, b_i \in \mathbb{R}.
\end{equation}
\noindent (iii) A connection between the generalized Wright function and the H-function is (see \cite{kilbas2002generalized}) 
\begin{equation}\label{gwf1}
_p\psi_q \left[z\; \vline \;\begin{matrix}
\left(\alpha_i, \beta_i\right)_{1,p}\\
(a_j,b_j)_{1,q}
\end{matrix} \right] =\;  H_{p,q+1}^{1,p} \left[-z\; \vline \;\begin{matrix}
\left(1-\alpha_i, \beta_i\right)_{1,p}\\
(0,1) & (a_j,b_j)_{1,q}
\end{matrix} \right],
\end{equation}
where the H-function $H_{p,q+1}^{1,p}(z)$  can be expressed in terms of Mellin-Barnes type integral (see \cite{kilbas2002generalized}). For $k \in \mathbb{N},$ the $k$th order partial derivative of the H-function is given by (see \cite{mathai2009h})
\begin{align}\label{deri1}
\diffp{^k}{z^k}&\left[z^{-la_1/b_1} H_{p_1,q_1}^{p_2,q_2} \left[z^l\; \vline \;\begin{matrix}
\left(\alpha_i, \beta_i\right)_{1,p_1}\\
(a_j,b_j)_{1,q_1}
\end{matrix} \right]\right] 
= z^{-k-(-la_1/b_1)}\left(\frac{-l}{b_1}\right)^k  H_{p_1,q_1}^{p_2,q_2} \left[z^l\; \vline \;\begin{matrix}
\left(\alpha_i, \beta_i\right)_{1,p_1}\\
(k+a_1, b_1), & (a_j,b_j)_{2,q_2}
\end{matrix} \right].
\end{align}
\subsection{Definitions and some elementary distributions}
\noindent (i)	Let $f :[a,b] \subset \mathbb{R}\longrightarrow\mathbb{R}$ be such that $f(t)$ is $(n+1)$ times continuous differentiable for $n < \tau <n+1$. Then, the Riemann-Liouville fractional derivative of order $\tau>0$ is defined as (see \cite{podlubny1999introduction})
\begin{equation*}
_aD^\tau_tf(t)=\bigg( \frac{d}{dt}\bigg)^{n+1}\int_{a}^{t}(t-u)^{n-\tau}f(u)du.
\end{equation*}
\noindent (ii) For  $ 0 < \beta \leq 1$, let  $\{N_\beta(t, \lambda)\}_{t\geq 0}$ be a FPP with parameter $\lambda>0$. Its one- dimensional distributions are given by (see \cite{laskin2003fractional, meerschaert2011fractional})
\begin{equation*}
p_{\beta}(n/t, \lambda) =P[N_\beta(t, \lambda)=n] =\frac{(\lambda t^{\beta})^n}{n!} \sum_{k=0}^{\infty} \frac{(n+k)!}{k!}\frac{(-\lambda t^{\beta})^k}{\Gamma(\beta(k+n)+1)},\;\; n\in \mathbb{Z}_{+}.
\end{equation*}
\noindent (iii) Let $\Gamma(t)\sim G(\mu, \rho t)$. Its probability density function (pdf) is given by 
\begin{equation*}
f_G(x,t) = \frac{\mu^{\rho t}}{\Gamma(\rho t)}x^{\rho t-1}e^{-\mu x}, \;\; x>0.
\end{equation*}
\noindent (iv) \cite{barndorff2000probability} and \cite{kumar2019tempered} discussed the Mittag-Leffler (ML) L\'{e}vy process with various properties.  For $\alpha \in (0,1)$ and $\rho , \mu, t > 0,$ let $\{M_{\alpha,\rho, \mu}(t)\}_{t \geq 0}$ be a ML L\'{e}vy process with L\'{e}vy measure density
\begin{equation}\label{lv1}
\pi (x) = \frac{\alpha \rho}{x} L_{\alpha,1}^1 \left(-\mu x^\alpha\right), \;\;x >0.
\end{equation}
The Laplace transform of  ML L\'{e}vy process is
\begin{equation}\label{lt1111}
\mathbb{E}\left(e^{-uM_{\alpha,\rho, \mu}(t)}\right) = \left(\frac{\mu}{\mu + u^\alpha}\right)^{\rho t}.
\end{equation}
Using the conditioning argument, the pdf of ML L\'{e}vy process is obtained and is given as (see \cite{ kumar2019tempered})
\begin{equation}\label{pdf1}
f_{M_{\alpha ,\rho, \mu}(t)}(x) = \sum_{k=0}^{\infty}(-1)^k\frac{\mu^{\rho t+k} \Gamma(\rho t+k)}{\Gamma(\rho t) \Gamma(k+1)}\frac{x^{\alpha(\rho t+k)-1}}{\Gamma(\alpha(\rho t+k))}, \; x> 0,
\end{equation}
The moments of $l$th order are of the following form
\begin{equation}\label{c1}
\mathbb{E}[M_{\alpha,\rho, \mu}^{l}(t)] = \frac{\rho t}{\mu^{l/\alpha} \Gamma(1-l)}B\left(1-\frac{l}{\alpha}, \rho t+\frac{l}{\alpha}\right) \sim \frac{\Gamma(1-\frac{l}{\alpha})}{\mu^{l/\alpha}\Gamma(1-l)}(\rho t)^{l/\alpha}, \; 0 < l<\alpha\; \text{as}\; t \rightarrow \infty,
\end{equation}
 where $B(z_1,z_2) = \frac{\Gamma(z_1)\Gamma(z_2)}{\Gamma(z_1 + z_2)},  \text{Re}(z_1) >0,  \text{Re}(z_2)>0 $ represents the beta function.

\section{Generalized Fractional Negative Binomial Process}
 Here, we define a generalized fractional negative binomial process $\{\mathcal{G}_{\alpha,\rho}^{\mu, \beta}(t,  \lambda)\}_{t\geq 0}$ by replacing the gamma subordinator with an independent ML L\'{e}vy subordinator in the gamma subordinated form of the FNBP, that is,
\begin{equation*}
\mathcal{G}_{\alpha,\rho}^{\mu, \beta}(t,  \lambda) := N_{\beta}(M_{\alpha,\rho, \mu}(t), \lambda), \;\; t\geq 0.
\end{equation*}
The probability mass function (pmf) of $\{\mathcal{G}_{\alpha,\rho}^{\mu, \beta}(t,  \lambda)\}_{t\geq 0}$, denoted by $p_{\alpha,\rho}^{\mu, \beta}(n,t)=P[\mathcal{G}_{\alpha,\rho}^{\mu, \beta}(t,  \lambda)=n]$ is derived as 
\begin{align*}
p_{\alpha,\rho}^{\mu, \beta}(n,t)
&= \int_{0}^{\infty} 	p_{\beta}(n/y, \lambda) 	f_{M_{\alpha,b, \mu}(t)}(y) dy\\
&= \int_{0}^{\infty} \left(\frac{(\lambda y^{\beta})^n}{n!} \sum_{k=0}^{\infty} \frac{(n+k)!}{k!}\frac{(-\lambda y^{\beta})^k}{\Gamma(\beta(k+n)+1)}\right)	f_{M_{\alpha,\rho, \mu}(t)}(y) dy\\
&= \frac{\lambda^n}{n!} \sum_{k=0}^{\infty} \frac{(n+k)!}{k!}\frac{(-\lambda )^k}{\Gamma(\beta(k+n)+1)} \int_{0}^{\infty} y^{\beta(n+k)}	f_{M_{\alpha,\rho, \mu}(t)}(y) dy\\
&= \frac{\lambda^n}{n!} \sum_{k=0}^{\infty} \frac{(n+k)!}{k!}\frac{(-\lambda )^k}{\Gamma(\beta(k+n)+1)} \mathbb{E}M_{\alpha,\rho, \mu}^{\beta(k+n)}(t),\;\;0 < \beta(k+n) <\alpha.
\end{align*}
\begin{remark}
		\noindent For $\alpha =1$, the pmf of GFNBP reduces to  
		\begin{equation*}
		p_{1,\rho}^{\mu, \beta}(n,t) =  \frac{\lambda^n}{n!} \sum_{k=0}^{\infty} \frac{(n+k)!}{\mu^{\beta(k+n)}k!}\frac{\Gamma((n+k)\beta+\rho t)(-\lambda )^k}{\Gamma(\rho t)\Gamma(\beta(k+n)+1)},
		\end{equation*}
		which is the pmf of the FNBP discussed in \cite{samy2018fractional}.
	\end{remark}
\begin{remark}
With the help of (\ref{c1}), the pmf $p_{\alpha,\rho}^{\mu, \beta}(n,t)$ can be expressed in terms of generalized Wright function of the following form
\begin{equation*}\label{pmf111}
p_{\alpha,\rho}^{\mu, \beta}(n,t) =\frac{\lambda^n}{\mu^{\beta n/\alpha} n! \Gamma(\rho t)}\;  _3\psi_2 \left[\frac{-\lambda}{\mu^{\beta/\alpha}}\; \vline \;\begin{matrix}
\left(n+1,1\right),& \left(1-\frac{\beta n}{\alpha}, \frac{-\beta}{\alpha}\right),& \left(\rho t + \frac{\beta n}{\alpha}, \frac{\beta}{\alpha}\right)\\
(1-\beta n, -\beta), &(1+\beta n, \beta)
\end{matrix} \right].
\end{equation*}
Alternatively, it may be also expressed in terms of the $H$-function via the relation 
\begin{equation}\label{hf45}
p_{\alpha,\rho}^{\mu, \beta}(n,t) = \frac{\lambda^n}{\mu^{\beta n/\alpha} n! \Gamma(\rho t)}H_{2,3}^{1,2} \left[\frac{\lambda}{\mu^{\beta/\alpha}}\; \vline \;\begin{matrix}
\left(-n,1\right),& \left(\frac{\beta n}{\alpha}, \frac{-\beta}{\alpha}\right),& \left(1- \rho t - \frac{\beta n}{\alpha}, \frac{\beta}{\alpha}\right)\\
(0,1),& (\beta n, -\beta), &(-\beta n, \beta)
\end{matrix} \right].
\end{equation}
\end{remark}
Applying the Leibniz rule of derivative for the convolution of functions in (\ref{hf45}) and with the help of (\ref{deri1}), the partial differential equations (PDEs) governed by  the pmf of GFNBP is
\begin{align*}
&\diffp{^k}{\mu^k} p_{\alpha,\rho}^{\mu, \beta}(n,t) =  \frac{1}{\mu^{\beta n/\alpha} n! \Gamma(bt)}\sum_{i =0}^{k}\binom{k}{i} \binom{n}{i}(-1)^{n-k} H_{2,3}^{1,2} \left[\frac{\lambda}{\mu^{\beta/\alpha}}\; \vline \;\begin{matrix}
\left(-n,1\right),& \left(\frac{\beta n}{\alpha}, \frac{-\beta}{\alpha}\right),& \left(1- \rho t - \frac{\beta n}{\alpha}, \frac{\beta}{\alpha}\right)\\
(k-i,1),& (\beta n, -\beta), &(-\beta n, \beta)
\end{matrix} \right],
\end{align*}
with
\begin{equation*}
p_{\alpha,\rho}^{\mu, \beta}(n,t)=%
\begin{cases}
1 &\text{if}\;\; n,t =0, \\
0 &\text{if}\;\; n \geq 1,\; t=0,
\end{cases} \text{ and } 	p_{\alpha,\rho}^{\mu, \beta}(n,t) = 0 \text{ for all }\; n< 0.
\end{equation*}
Next, we show that the pmf of the  GFNBP satisfies a fractional PDE. 
\begin{lemma}\label{lemm1}
	(\cite{samy2018fractional}) For $0 <\tau \leq 1$, the governing fractional PDE for the gamma subordinator $\{\Gamma(t)\}_{t\geq 0}$ is given by
	\begin{equation*}
	\diffp{^\tau}{t^\tau}	f_G(x,t) = \rho	\diffp{^{\tau-1}}{t^\tau}\left[\log \alpha +\log y - \psi(\rho 
	t)\right]	f_G(x,t), \;\; y >0 \text{ and }
	f_G(x,0) =0,
	\end{equation*}
	where $\psi(x)$ is the digamma function and $	\diffp{^\tau}{t^\tau}(\cdot)$ is the Riemann-Liouville fractional differential operator.
\end{lemma}
The next lemma gives the fractional version of PDE with respect to time variable satisfying the pdf of the ML L\'{e}vy process.
\begin{lemma}\label{prop1}
	Let $	g_\alpha(x,t)$ be the pdf for the $\alpha$-stable process. 
	Then, the density of the ML L\'{e}vy process satisfies the following fractional PDE
	\begin{equation*}
	\diffp{^\tau}{t^\tau}  f_{M_{\alpha,\rho, \mu}(t)}(x) = \rho	\diffp{^{\tau-1}}{t^{\tau-1}}\left[\left(\log \alpha - \psi(\rho t) \right) f_{M_{\alpha,\rho, \mu}(t)}(x)+\int_{0}^{\infty}g_\alpha(x,y) (\log y) f_{M_{\alpha,\rho, \mu}(t)}(y)dy\right].
	\end{equation*}
\end{lemma}
\begin{proof} Consider
	\begin{equation*}
	f_{M_{\alpha,\rho, \mu}(t)}(x) 
	= \int_{0}^{\infty} g_\alpha(x,y)	f_G(y,t) dy.
	\end{equation*}
Operating the Riemann-Liouville fractional derivative, we get
	\begin{align*}
	\diffp{^\tau}{t^\tau} f_{M_{\alpha,\rho, \mu}(t)}(x) 
	&= \diffp{^\tau}{t^\tau} \int_{0}^{\infty} g_\alpha(x,y)	f_G(y,t) dy\\
	&= \int_{0}^{\infty} g_\alpha(x,y) 	\diffp{^\tau}{t^\tau} 	f_G(y,t) dy\\
	&= \int_{0}^{\infty} g_\alpha(x,y) \left[\rho \diffp{^{\tau-1}}{t^\tau}\left[\log \alpha +\log y - \psi(\rho t)\right]	f_G(y,t)\right]dy \;\;\;(\text{using  Lemma \ref{lemm1}})\\
	&= \rho	\diffp{^{\tau-1}}{t^{\tau-1}}\int_{0}^{\infty} g_\alpha(x,y)\left(\log \alpha - \psi(\rho t) \right) f_G(y,t)dy+\rho \int_{0}^{\infty}g_\alpha(x,y) (\log y) 	\diffp{^{\tau-1}}{t^{\tau-1}}f_G(y,t)dy.
	\end{align*}
With the help of simple algebra, the lemma follows.
\end{proof}
Now, using Lemma  \ref{prop1}, we can get the governing fractional PDE for the GFNBP with respect to time variable of the form
\begin{align*}
&\frac{1}{\rho}	\diffp{^\tau}{t^\tau} 	p_{\alpha,\rho}^{\mu, \beta}(n,t) = 	\diffp{^{\tau-1}}{t^{\tau-1}}\left[\left(\log \alpha - \psi(\rho t) \right) 	p_{\alpha,\rho}^{\mu, \beta}(n,t)+\int_{0}^{\infty}\int_{0}^{\infty}	p_{\beta}(n/s, \lambda)g_\alpha(x,y) (\log y) f_{M_{\alpha,\rho, \mu}(t)}(y)dy ds\right],\;\text{with}\; p_{\alpha,\rho}^{\mu, \beta}(0,0) = 1.
\end{align*}
\subsection{Mean, variance, autocovariance and index of dispersion}
\begin{theorem}\label{thm1} Let  $\{\mathcal{G}_{\alpha,\rho}^{\mu, \beta}(t,  \lambda)\}_{t\geq 0}$ be a GFNBP. For $0 <s \leq t <\infty$, we have \\
	\noindent (i)  $\mathbb{E}[\mathcal{G}_{\alpha,\rho}^{\mu, \beta}(t,  \lambda)] = q\mathbb{E}[M_{\alpha,\rho, \mu}^{\beta}(t)] \sim \frac{q \Gamma(1-\frac{\beta}{\alpha})}{\mu^{\beta/\alpha}\Gamma(1-\beta)}(\rho t)^{\beta/\alpha}, \; 0<\beta < \alpha.$\\
\noindent	(ii) $\text{Var}[\mathcal{G}_{\alpha,\rho}^{\mu, \beta}(t,  \lambda)] =  q\mathbb{E}[M_{\alpha,\rho, \mu}^{\beta}(t)] - q^2\left(\mathbb{E}[M_{\alpha,\rho, \mu}^{\beta}(t)]\right)^2 + 2d  \mathbb{E}[M_{\alpha,\rho, \mu}^{2\beta}(t)].$\\
		\noindent (iii) Cov$[\mathcal{G}_{\alpha,\rho}^{\mu, \beta}(s,  \lambda), \mathcal{G}_{\alpha,\rho}^{\mu, \beta}(t,  \lambda)]  = q\mathbb{E}[M_{\alpha,\rho, \mu}^{\beta}(t)] +d\mathbb{E}[M_{\alpha,\rho, \mu}^{2 \beta}(s)]-q^2 \mathbb{E}[M_{\alpha,\rho, \mu}^{\beta}(s)]\mathbb{E}[M_{\alpha,\rho, \mu}^{\beta}(t)]+ q^2\beta \mathbb{E}\left[M_{\alpha,\rho, \mu}^{2\beta}(t)B\left(\beta, 1+\beta; \frac{M_{\alpha,\rho, \mu}(s)}{M_{\alpha,\rho, \mu}(t)} \right)\right],$
	where $q = \frac{\lambda}{\Gamma(1+\beta)}, d = \beta q^2 B(\beta, 1+\beta),$ and $ B(r,s;x) = \int_{0}^{x} t^{r-1} (1-t)^{s-1}dt$ for $ 0 < x< 1$ is an incomplete beta function.
\end{theorem}
\begin{proof}
	The mean, variance, and autocovariance functions of the FPP is given by (see \cite{laskin2003fractional})
	\begin{equation}\label{ev1}
	\mathbb{E}\left[N_\beta(t, \lambda)\right] = q t^\beta, \;\;\;\;\text{Var}\left[N_\beta(t, \lambda)\right]=qt^\beta + Rt^{2\beta},
	\end{equation}
	\begin{equation*}
	\text{Cov}\left[N_\beta(s, \lambda), N_\beta(t, \lambda)\right] = qs^\beta +ds^{2\beta}+q^2\left[\beta t^{2\beta} B(\beta, 1+\beta;s/t)- (st)^{\beta}\right],\;\; 0< s \leq t.
	\end{equation*}
	Using the conditioning argument and with the help of (\ref{ev1}) and (\ref{c1}), Part (i) of the theorem can be easily  obtained.
	Also, one may derive that
	\begin{equation*}
	\mathbb{E}\left[N_\beta(s, \lambda)N_\beta(t, \lambda)\right] = qs^\beta +ds^{2\beta}+ q^2\left[\beta t^{2\beta} B(\beta, 1+\beta;s/t)\right],
	\end{equation*}
	that gives
	\begin{align*}\label{cov1}
	\mathbb{E}\left[\mathcal{G}_{\alpha,\rho}^{\mu, \beta}(s,  \lambda) \mathcal{G}_{\alpha,\rho}^{\mu, \beta}(t,  \lambda)\right] 
	&= \mathbb{E}\left[\mathbb{E}\left[N_\beta(M_{\alpha,\rho, \mu}(s))N_\beta(M_{\alpha,\rho, \mu}(t))|(M_{\alpha,\rho, \mu}(s), M_{\alpha,\rho, \mu}(t))\right]\right] \nonumber\\
	&=  q\mathbb{E}[M_{\alpha,\rho, \mu}^{\beta}(t)] +d\mathbb{E}[M_{\alpha,\rho, \mu}^{2 \beta}(s)]+ q^2\beta \mathbb{E}\left[M_{\alpha,\rho, \mu}^{2\beta}(t)B\left(\beta, 1+\beta; \frac{M_{\alpha,\rho, \mu}(s)}{M_{\alpha,\rho, \mu}(t)} \right)\right].
	\end{align*}
Now, the covariance formula gives  Part (iii) of the theorem. Part (ii) is a consequence of Part (iii), when $s=t$.
\end{proof}
To study the index of dispersion of GFNBP, we first prove the following lemma.
\begin{lemma}\label{l1}
	For $0 < l < \alpha/2,$ we have
	\begin{equation*}
	\left(\mathbb{E}\left[ M_{\alpha,b, \mu}^{l}(t) \right]\right)^2 \leq 	\mathbb{E}\left[M_{\alpha,b, \mu}^{2l}(t) \right].
	\end{equation*}
\end{lemma}
\begin{proof} By self-similar property of the stable processes, we get
	\begin{equation*}
	\mathbb{E}[M_{\alpha,\rho, \mu}^{2l}(t)] = \mathbb{E}[G_{\mu,\rho}^{2l/\alpha}(t)]\mathbb{E}[S_\alpha^{2l}(1)] \geq 	\mathbb{E}[G_{\mu,\rho}^{l/\alpha}(t)]^2 	\mathbb{E} [S_\alpha^{l}(1)]^2 = \left(\mathbb{E}[M_{\alpha,\rho, \mu}^{l}(t)]\right)^2.
	\end{equation*}
\end{proof}
\noindent A stochastic process $\{X(t)\}_{t \geq 0}$ is overdispesed if $Var[X(t)]-\mathbb{E}[X(t)]>0$ for all $t\geq 0$ (see \cite[p. 72]{cox1966statistical}).
Using Lemma \ref{l1}, we have
\begin{align*}
\text{Var}[\mathcal{G}_{\alpha,\rho}^{\mu, \beta}(t,  \lambda)] - \mathbb{E}[\mathcal{G}_{\alpha,\rho}^{\mu, \beta}(t,  \lambda)]  
&= 2d \mathbb{E}[M_{\alpha,\rho, \mu}^{2\beta}(t)]- \left(q\mathbb{E}[M_{\alpha,\rho, \mu}^{\beta}(t)]\right)^2
= \frac{\lambda^2}{\beta} \left(\frac{\mathbb{E}\left[M_{\alpha,\rho, \mu}^{2\beta}(t)\right]}{\Gamma(2 \beta)} - \frac{\left(\mathbb{E}\left[M_{\alpha,\rho, \mu}^{\beta}(t)\right]\right)^2}{\beta \Gamma^2 (\beta)}\right) \geq 0.
\end{align*}
It is clear from the fact that $\frac{\lambda^2}{\beta} \left(\frac{1}{\Gamma(2 \beta)} - \frac{1}{\beta \Gamma^2 (\beta)}\right) > 0$ for $\lambda >0$ and $\beta \in (0,1)$ (see \cite{beghin2014fractional}). Therefore, GFNBP exhibits overdispersion.
\subsection{Laplace transform}Let $q(x,t)$ be the pdf of the $E_\beta\left(M_{\alpha,\rho, \mu}(t)\right)$ and $h_\beta(x,t)$ be the pdf of the inverse stable subordinator $E_\beta(t)$ with Laplace transform $\mathbb{E}[e^{-uE_\beta(t)}] =L_{\beta,1}^{1}(-ut^\beta)$ (see \cite{meerschaert2013inverse}). Then, the Laplace transform of    $E_\beta\left(M_{\alpha,\rho, \mu}(t)\right)$ can be derived as
\begin{align*}
\mathbb{E}\left[e^{-u E_\beta\left(M_{\alpha,\rho, \mu}(t)\right)}\right] &= \int_{0}^{\infty} e^{-ux}q(x,t) = \int_{0}^{\infty} \int_{0}^{\infty} e^{-ux} h_\beta(x,y) f_{M_{\alpha,\rho, \mu}(t)}(y) dy dx\\
&= \int_{0}^{\infty} L_{\beta,1}^{1}(-uy^\beta) f_{M_{\alpha,\rho, \mu}(t)}(y) dy\\
&= \sum_{k=0}^{\infty}  \frac{(-u)^k}{\Gamma(1+\beta k)} \mathbb{E}[M_{\alpha,\rho, \mu}^{k\beta}(t)]\;\; (0 <\beta k < \alpha)\\
&= \frac{1}{\Gamma(\rho t)}\;  _2\psi_2 \left[\frac{-u}{\mu^{\beta/\alpha}}\; \vline \;\begin{matrix}
\left(1, -\frac{\beta}{\alpha}\right), & \left(\rho t, \frac{\beta }{\alpha}\right)\\
(1,\beta), & (1, -\beta)
\end{matrix} \right].
\end{align*}
Using the conditioning arguments,  we obtain the Laplace transform for the GFNBP as
\begin{align*}
\mathbb{E}\left[e^{-u\mathcal{G}_{\alpha,\rho}^{\mu, \beta}(t)}\right] &= \mathbb{E}\left[\mathbb{E}\left[\text{exp}\left(-\lambda E_\beta\left( M_{\alpha,\rho, \mu}(t)\right)(1-e^{-u})\right)/  E_\beta\left( M_{\alpha,\rho, \mu}(t)\right)\right]\right]\\
&= \frac{1}{\Gamma(\rho t)}\;  _2\psi_2 \left[\frac{-\lambda (1-e^{-u}) }{\mu^{\beta/\alpha}}\; \vline \;\begin{matrix}
\left(1, -\frac{\beta}{\alpha}\right), & \left(\rho t, \frac{\beta }{\alpha}\right)\\
(1,\beta), & (1, -\beta)
\end{matrix} \right].
\end{align*}
\begin{remark}
The probability generating function of the GFNBP can be deduced from the Laplace transform and is given by
\begin{equation*}
\mathbb{E}\left[u^{\mathcal{G}_{\alpha,\rho}^{\mu, \beta}(t)}\right] = \frac{1}{ \Gamma(\rho t)}\;  _2\psi_2 \left[\frac{-\lambda(1-u)}{\mu^{\beta/\alpha}}\; \vline \;\begin{matrix}
\left(1, \frac{-\beta}{\alpha}\right),& \left(\rho t, \frac{\beta}{\alpha}\right)\\
(1, \beta), &(1, -\beta)
\end{matrix} \right].
\end{equation*} 
\end{remark}
\subsection{Infinite divisibility}
The following self-similarity properties of stable and inverse stable subordinators will be used (see \cite{meerschaert2004limit})

\begin{equation}\label{de1}
S_\alpha(t)  \stackrel{d}{=} t^{1/\alpha}S_\alpha(1) \text{ and } E_{\beta}(t) \stackrel{d}{=}t^{\beta} E_{\beta}(1),
\end{equation}
where $\stackrel{d}{=}$ stands for equality in distribution.
Also, we may observe that
\begin{equation*}
N_{\beta}(t, \lambda) \stackrel{d}{=} N\left(E_{\beta}(t), \lambda\right)\stackrel{d}{=} N\left(t^{\beta}E_{\beta}(1), \lambda\right).
\end{equation*}
By the renewal theorem (see \cite{samy2018fractional}), we get
\begin{equation*}\label{e1}
\lim_{t \rightarrow \infty} \frac{N\left(t^{\beta}E_{\beta}(1), \lambda\right)}{t^{\beta}} = \lambda E_{\beta}(1)\;a.s.
\end{equation*}

\begin{theorem}\label{thm2}
	The GFNBP $\{\mathcal{G}_{\alpha,\rho}^{\mu, \beta}(t,  \lambda)\}_{t\geq 0}$ is not infinitely divisible.
\end{theorem}
\begin{proof}
	Using (\ref{de1}), we get
	\begin{align*}
	\mathcal{G}_{\alpha,\rho}^{\mu, \beta}(t,  \lambda) &= N_{\beta}(M_{\alpha,\rho, \mu}(t), \lambda) 
	= N_{\beta}(S_\alpha(G_{\mu,\rho}(t)), \lambda)
	= N\left(\left[S_\alpha G_{\mu,\rho}(t)\right]^{\beta}E_{\beta}(1), \lambda\right)\\
	&= N\left(\left\{\left[G_{\mu,\rho}(t)\right]^{1/\alpha} S_\alpha(1)\right\}^{\beta}E_{\beta}(1), \lambda\right)
	= N\left(\left[G_{\mu,\rho}(t)\right]^{\beta/\alpha}\left( S_\alpha(1)\right)^{\beta}E_{\beta}(1), \lambda\right).
	\end{align*}
	Now, we consider
	\begin{align*}
~~~~~~	& \lim_{t \rightarrow \infty} \frac{N\left(\left[G_{\mu,\rho}(t)\right]^{\beta/\alpha}\left( S_\alpha(1)\right)^{\beta}E_{\beta}(1), \lambda\right)}{t^{\beta/\alpha}}\\
	&\stackrel{d}{=}  \lim_{t \rightarrow \infty} \frac{N\left(\left[G_{\mu,\rho}(t)\right]^{\beta/\alpha}\left( S_\alpha(1)\right)^{\beta}E_{\beta}(1), \lambda\right)}{\left[G_{\mu,\rho}(t)\right]^{\beta/\alpha}\left( S_\alpha(1)\right)^{\beta}E_{\beta}(1)} \frac{\left[G_{\mu,\rho}(t)\right]^{\beta/\alpha}\left( S_\alpha(1)\right)^{\beta}E_{\beta}(1)}{t^{\beta/\alpha}}\\
	&\stackrel{d}{=}  \;\;\; \lambda E_{\beta}(1)\left( S_\alpha(1)\right)^{\beta}  \lim_{t \rightarrow \infty} \left(\frac{G_{\mu,\rho}(t)}{t}\right)^{\beta/\alpha}\\
	&\stackrel{d}{=}  \;\;\; \lambda E_{\beta}(1)\left( S_\alpha(1)\right)^{\beta}  \left(\mathbb{E}G_{\mu,\rho}(1)\right)^{\beta/\alpha}.
	\end{align*}
	For a large t,  $G_{\mu,\rho}(t)/{t} \rightarrow \mathbb{E}G_{\mu,\rho}(1).$  Since $E_{\beta}(1)$ is not infinitely divisible, the result follows.
\end{proof}
\subsection{Dependence structure}
\begin{definition}
	For $0 <s <t,$ let the correlation function Corr$[X(s), X(t)]$ for a stochastic process $\{X(t)\}_{t \geq 0}$ satisfies the following relation
	\begin{equation*}
	c_1(s)t^{-d} \leq \text{ Corr}[X(s), X(t)]\leq c_2(s)t^{-d}
	\end{equation*}
	for large  $t$, $d >0$, $c_1(s) > 0$ and $c_2(s) > 0.$ Expressly
	\begin{equation*}
	\lim_{t \rightarrow \infty} \frac{ \text{ Corr}[X(s), X(t)]}{t^{-d}} = c(s),
	\end{equation*}
	for some $c(s) >0$ and $d > 0$.
	The process $\{X(t)\}_{t \geq 0}$ is said to have the LRD property if $d \in (0,1)$.
\end{definition}
\noindent The following lemma can be proved exactly in a similar fashion as Lemma 2 in \cite{maheshwari2016long}.
\begin{lemma}\label{lm22}
	Let $\beta \in (0,\alpha)$ and $0 < s < t, s$ is fixed. Then,  the following asymptotic expansion holds for a large t.
	\begin{align*}
	&\text{(i)}\;\;\mathbb{E}\left[M_{\alpha,\rho, \mu}^{\beta}(s) M_{\alpha,\rho, \mu}^{\beta}(t) \right] \sim \mathbb{E}\left[M_{\alpha,\rho, \mu}^{\beta}(s)\right] \mathbb{E}\left[M_{\alpha,\rho, \mu}^{\beta}(t-s)\right].\\
	&\text{(ii)}\;\;\beta \mathbb{E}\left[M_{\alpha,\rho, \mu}^{2\beta}(t)B\left(\beta, 1+\beta; \frac{M_{\alpha,\rho, \mu}(s)}{M_{\alpha,\rho, \mu}(t)} \right)\right] \sim \mathbb{E}\left[M_{\alpha,\rho, \mu}^{\beta}(s)\right] \mathbb{E}\left[M_{\alpha,\rho, \mu}^{\beta}(t-s)\right].
	\end{align*}
\end{lemma}

\begin{theorem}
	The GFNBP exhibits the LRD property.
\end{theorem}
\begin{proof} Using (\ref{c1}) and with the help of Lemma \ref{lm22}(ii), the asymptotic behaviour of Theorem \ref{thm1}(iii) is 
	\begin{align*}
	 \text{Cov}[\mathcal{G}_{\alpha,\rho}^{\mu, \beta}(s,  \lambda), \mathcal{G}_{\alpha,\rho}^{\mu, \beta}(t,  \lambda)] &\sim q\mathbb{E}[M_{\alpha,\rho, \mu}^{\beta}(t)] +d\mathbb{E}[M_{\alpha,\rho, \mu}^{2 \beta}(s)]-q^2 \mathbb{E}[M_{\alpha,\rho, \mu}^{\beta}(s)]\left[\mathbb{E}[M_{\alpha,\rho, \mu}^{\beta}(t)] - \mathbb{E}\left[M_{\alpha,\rho, \mu}^{\beta}(t-s)\right]\right]\\
	 &\sim  q\mathbb{E}[M_{\alpha,\rho, \mu}^{\beta}(t)] +d\mathbb{E}[M_{\alpha,\rho, \mu}^{2 \beta}(s)]-q^2 \mathbb{E}[M_{\alpha,\rho, \mu}^{\beta}(s)]\left[\frac{\Gamma(1-\frac{\beta}{\alpha})}{\mu^{\beta/\alpha}\Gamma(1-\beta)}(\rho t)^{\beta/\alpha} - \frac{\Gamma(1-\frac{\beta}{\alpha})}{\mu^{\beta/\alpha}\Gamma(1-\beta)}(\rho (t-s))^{\beta/\alpha}\right]\\
	 &\sim  q\mathbb{E}[M_{\alpha,\rho, \mu}^{\beta}(t)] +d\mathbb{E}[M_{\alpha,\rho, \mu}^{2 \beta}(s)]. \;\;\; (\text{since}\;\; t^{\beta/\alpha} - (t-s)^{\beta/\alpha} \sim \beta \alpha^{-1} st^{(\beta/\alpha) -1}) 
		\end{align*}
		Also, the asymptotic behaviour of Theorem \ref{thm1}(ii) follows
		\begin{align*}
		\text{Var}\left[ \mathcal{G}_{\alpha,\rho}^{\mu, \beta}(t, \lambda)\right]
		&\sim \left(\frac{\rho t}{\mu}\right)^{2\beta / \alpha}\left(\frac{2d\Gamma(1-\frac{2\beta}{\alpha})}{\Gamma(1-2\beta)} - \left(\frac{q\Gamma(1-\frac{\beta}{\alpha})}{\Gamma(1-\beta)}\right)^2 \right) \sim t^{2\beta/\alpha}d_1,
		\end{align*}
where $d_1 = 	\left(\frac{\rho}{\mu}\right)^{2\beta / \alpha}\left(\frac{2d\Gamma(1-\frac{2\beta}{\alpha})}{\Gamma(1-2\beta)} - \left(\frac{q\Gamma(1-\frac{\beta}{\alpha})}{\Gamma(1-\beta)}\right)^2 \right)$.
Therefore, we have the correlation function as
	\begin{align*}
	\text{Corr}\left[ \mathcal{G}_{\alpha,\rho}^{\mu, \beta}(s),  \mathcal{G}_{\alpha,\rho}^{\mu, \beta}(t)\right] &= \frac{\text{Cov}\left[ \mathcal{G}_{\alpha,\rho}^{\mu, \beta}(s),  \mathcal{G}_{\alpha,\rho}^{\mu, \beta}(t)\right]}{\sqrt{\text{Var}\left[ \mathcal{G}_{\alpha,\rho}^{\mu, \beta}(s)\right]}\sqrt{\text{Var}\left[ \mathcal{G}_{\alpha,\rho}^{\mu, \beta}(t)\right]}} \sim t^{-\beta/\alpha}\left(\frac{q\mathbb{E}[M_{\alpha,\rho, \mu}^{\beta}(s)] + d\mathbb{E}[M_{\alpha,\rho, \mu}^{2\beta}(s)]}{\sqrt{d_1\text{Var}\left[ \mathcal{G}_{\alpha,\rho}^{\mu, \beta}(s)\right]}}\right).
	\end{align*}
	Hence, for $0< \beta < \alpha$ and the decaying power $t^{-\beta/\alpha}$,  the GFNBP has the LRD property.
\end{proof}
\section{Space Fractional and Non-Homogeneous Versions}
\subsection{Space fractional version of the GFNBP}
\noindent  \cite{orsingher2012space} have studied the space fractional Poisson process (SFPP) which is characterized as the Poisson process time-changed by independent stable subordinator. Let $\{N_{\alpha '}(t, \lambda)\}_{t \geq 0}, \; 0< \alpha ' <1$ be a SFPP with pmf 
\begin{equation}\label{stfpp1}
p_{\alpha'}(n/t, \lambda) = \frac{(-1)^n}{n!} \sum_{k=0}^{\infty} \frac{(- \lambda^{\alpha'} t)^k}{k!} \frac{\Gamma(k{\alpha'}+1)}{\Gamma(k{\alpha'}+1-n)}, \; \lambda >0.
\end{equation}
Here,  we present a space fractional version of the GFNBP by subordinating the SFPP with an independent ML L\'{e}vy subordinator. We denote it by $\{\mathcal{H}_{\alpha,\rho}^{\mu, \alpha '}(t,  \lambda)\}_{t \geq 0}$ and is defined as
\begin{equation*}
\mathcal{H}_{\alpha,\rho}^{\mu, \alpha '}(t,  \lambda) = N_{\alpha '}(M_{\alpha,\rho, \mu}(t), \lambda), \;\; t\geq 0.
\end{equation*}  

Using (\ref{lv1}) and the formula given in p. 197 of \cite{ken1999levy}, one may calculate the L\'{e}vy measure density $\mathcal{V}$ for the process $\{\mathcal{H}_{\alpha,\rho}^{\mu, \alpha '}(t,  \lambda)\}_{t \geq 0}$  as
\begin{align}\label{lm56}
\mathcal{V}(k) 
&= \int_{0}^{\infty} \sum_{i=1}^{\infty} p_{\alpha'}(i/t, \lambda)  \delta_{\{i\}}(k)\pi(t) dt  \nonumber \\
&= \int_{0}^{\infty} \sum_{i=1}^{\infty} \left(  \frac{(-1)^i}{i!} \sum_{j=0}^{\infty} \frac{(- \lambda^{\alpha'} t)^j}{j!} \frac{\Gamma(j{\alpha'}+1)}{\Gamma(j{\alpha'}+1-i)}\right) \delta_{\{i\}}(k)\pi(t) dt  \nonumber\\
&= \rho \sum_{i=1}^{\infty}  \frac{(-1)^i}{i!}\delta_{\{i\}}(k) \sum_{j=0}^{\infty}  \frac{(- \lambda^{\alpha'})^j}{j!} \frac{\Gamma(j{\alpha'}+1)}{\Gamma(j{\alpha'}+1-i)} \int_{0}^{\infty} \alpha t^{j-1}L_{\alpha, 1}^1(-\mu t^\alpha) dt \nonumber \\
&= \rho \sum_{i=1}^{\infty}  \frac{(-1)^i}{i!}\delta_{\{i\}}(k)\sum_{j=0}^{\infty} \frac{(- \lambda^{\alpha'} /\mu ^{1/\alpha})^j}{j!} \frac{\Gamma(j{\alpha'}+1) \Gamma(j/\alpha)\Gamma(1-j/\alpha)}{ \Gamma(j{\alpha'}+1-i)\Gamma(1- j)}. 
\end{align}
The (\ref{lm56}) is obtained by an application of  the Mellin transform integral formula (see  \cite{shukla2007generalization}).
Also, in terms of the generalized Wright function, the L\'{e}vy measure density $\mathcal{V}$ can be re-written as
	\begin{equation*}
	\mathcal{V}(k) = \rho \sum_{i=1}^{\infty}  \frac{(-1)^i}{i!}\delta_{\{i\}}(k) \;  _3\psi_2 \left[\frac{-\lambda^{\alpha '}}{\mu^{1/\alpha}}\; \vline \;\begin{matrix}
	\left(1, \alpha '\right),& \left(0, \frac{1}{\alpha}\right),& \left(1,- \frac{1}{\alpha}\right)\\
	(1-i, \alpha), &(1, -1)
	\end{matrix} \right].
	\end{equation*}
	
\begin{remark}
	For $\alpha' =1,$ the L\'{e}vy measure $	\mathcal{V}$  coincides with the L\'{e}vy measure for the space fractional negative binomial process as reported in  \cite{beghin2018space}.
\end{remark}
It may be noted that the space fractional version of the GFNBP is again a subordinator. The Laplace transform for the SFPP  $\{N_{\alpha '}(t, \lambda)\}_{t \geq 0}$ is given by (see \cite{orsingher2012space})
\begin{equation*}
	\mathbb{E}\left[e^{-uN_{\alpha '}(t, \lambda) }\right] = \sum_{n=0}^{\infty}e^{-u n} p_{\alpha'}(n/t, \lambda) = e^{-t \lambda^{\alpha '}(1-e^{-u})^{\alpha '}}.
\end{equation*}
Using (\ref{lt1111}) and (\ref{stfpp1}), the Laplace transform for the distribution of the process  $\{\mathcal{H}_{\alpha,\rho}^{\mu, \alpha '}(t,  \lambda)\}_{t \geq 0}$ is calculated as
\begin{align*}
	\mathbb{E}\left[e^{-u\mathcal{H}_{\alpha,\rho}^{\mu, \alpha '}(t,  \lambda) }\right] &= \int_{0}^{\infty} f_G(y,t) \sum_{n=0}^{\infty} e^{-un} p_{\alpha'}(n/t, \lambda) dy 
	= \int_{0}^{\infty}  e^{-t \lambda^{\alpha '}(1-e^{-u})^{\alpha '}} f_G(y,t) dy\\
	&= \exp\left\{-\rho t \ln\left(1+ \frac{ \lambda^{\alpha '\alpha}(1-e^{-u})^{\alpha ' \alpha}}{\alpha} \right)\right\}.
\end{align*}
Hence, the Laplace exponent $\Upsilon_ {\mathcal{H}_{\alpha,\rho}^{\mu, \alpha '}(t,  \lambda)}$ of the process $\{\mathcal{H}_{\alpha,\rho}^{\mu, \alpha '}(t,  \lambda)\}_{t \geq 0}$  takes the following form
\begin{equation*}
\Upsilon_ {\mathcal{H}_{\alpha,\rho}^{\mu, \alpha '}(t,  \lambda)}(u) = -\frac{1}{t} \ln\left(	\mathbb{E}\left[e^{-u\mathcal{H}_{\alpha,\rho}^{\mu, \alpha '}(t,  \lambda) }\right]\right) = \rho \ln \left(1+ \frac{ \lambda^{\alpha '\alpha}(1-e^{-u})^{\alpha ' \alpha}}{\alpha}\right).
\end{equation*}
\subsection{Non-homogeneous version of the GFNBP} \cite{orsingher2012space} studied a space-time fractional Poisson process (STFPP) and discussed its connection to the associated fractional PDE. \cite{maheshwari2019non} introduced the non-homogeneous version of the STFPP. In this subsection, we present a non-homogeneous version of the GFNBP.

 Let  $\{	\hat{N}_{\alpha', \beta'}(t)\}_{t\geq0}$ be the non-homogeneous STFPP with rate function $\mathcal{R}(t) = \int_{0}^{t} \lambda(s)ds,$ where
$\lambda(s), \; s>0$ is intensity varying over time. The non-homogeneous version of the STFPP is defined as
\begin{equation}\label{nh11}
\hat{N}_{\alpha', \beta'}(\mathcal{R}(t), 1) \stackrel{d}{=} N\left(S_{\alpha'}(E_{\beta'}(\mathcal{R}(t))),1\right),
\end{equation} 
with pmf as
\begin{equation*}
\hat{p}_{\beta'}^{\alpha'}(n;\mathcal{R}(t)) = \frac{(-1)^n}{n!} \sum_{k=0}^{\infty} \frac{(- (\mathcal{R}(t) )^{\beta'})^k}{\Gamma(k {\beta'} +1)} \frac{\Gamma(k{\alpha'}+1)}{\Gamma(k{\alpha'}+1-n)}.
\end{equation*}
For $\beta'=1,$  (\ref{nh11}) coincides with the pmf of non-homogeneous version of the SFPP. Also,  when $\alpha'=1,$  (\ref{nh11}) corresponds to pmf of the non-homogeneous version of TFPP.

Now, we define a non-homogeneous STFPP time changed by the ML L\'{e}vy subordinator as
\begin{equation}\label{mm1}
\hat{\mathcal{G}} (t, \lambda) \stackrel{d}{=} 	\hat{N}_{\alpha', \beta'}(\mathcal{R}(M_{\alpha,\rho, \gamma}(t)), 1) \stackrel{d}{=} N\left(S_{\alpha'}(E_{\beta'}(\mathcal{R}(M_{\alpha,\rho, \mu}(t)))),1\right).
\end{equation}
The pmf for the process $\{\hat{\mathcal{G}} (t, \lambda)\}_{t\geq0}$ is computed as
\begin{align*}
\hat{P}(n;\mathcal{R}(t)) 
&= \int_{0}^{\infty} \hat{p}_{\beta'}^{\alpha'}(n;\mathcal{R}(y)) 	f_{M_{\alpha,\rho, \mu}(t)}(y) dy \\
&= \frac{(-1)^n}{n!} \sum_{k=0}^{\infty} \frac{(-1)^k}{\Gamma(k \beta' +1)} \frac{\Gamma(k\alpha'+1)}{\Gamma(k\alpha'+1-n)} \int_{0}^{\infty}  \mathcal{R}^{\beta'} (y)	f_{M_{\alpha,\rho, \mu}(t)}(y) dy\\
&= \frac{(-1)^n}{n!} \sum_{k=0}^{\infty} \frac{(-1)^k}{\Gamma(k \beta' +1)} \frac{\Gamma(k\alpha'+1)}{\Gamma(k\alpha'+1-n)} \mathbb{E}\left[\mathcal{R}^{\beta'} (M_{\alpha,\rho, \mu}(t))\right],\;\text{provided}\;\mathbb{E}\left[\mathcal{R}^{\beta'} (M_{\alpha,\rho, \mu}(t))\right] < \infty.
\end{align*}

\begin{remark}
	It is noted that, if we choose $\alpha=1,$ then pmf of process defined in (\ref{mm1}) coincides with the pmf of the non-homogeneous version of the space-time fractional negative binomial process (see \cite{maheshwari2019non}).
	In addition, when $\beta'=1,$ the pmf of (\ref{mm1}) leads to pmf of the non-homogeneous version of the space fractional negative binomial process (see \cite{beghin2018space}). Moreover, the  non-homogeneous version of the GFNBP can be considered as a  limiting case of (\ref{mm1}) when $\alpha'$ approaches to one.
\end{remark}
With the help of Lemma \ref{lemm1} and Lemma \ref{prop1}, we obtain the following formof  the fractional PDE satisfying the pmf of the non-homogeneous version of the space-time fractional negative binomial process (STFNBP) 
\begin{align*}
\diffp{^\tau}{t^\tau} 	\hat{p} (n;\mathcal{R}(t))  = \rho	\diffp{^{\tau-1}}{t^{\tau-1}}\left[\left(\log \alpha - \psi(\rho t) \right)	\hat{p} (n;\mathcal{R}(t))+\int_{0}^{\infty}	\hat{p}(n;\mathcal{R}(y))  (\log y) f_{M_{1,\rho, \mu}(t)}(y)dy\right].
\end{align*}
\begin{example} Let $\mathcal{R}(t) = \lambda^{\alpha'/{\beta'}}t.$ Using the self-similarity property of stable and inverse stable subordinators, we observe that, (\ref{mm1}) reduces to time changed version of the STFPP as
	\begin{align}\label{nhstfpp1}
	\hat{\mathcal{G}}(t, \lambda) 
	&\stackrel{d}{=} N\left(S_{\alpha'}(E_{\beta'}(\mathcal{R}(M_{\alpha,\rho, \mu}(t)))),1\right) \nonumber
	\stackrel{d}{=}  N\left(S_{\alpha'}(E_{\beta'}(\lambda^{\alpha'/\beta'}(M_{\alpha,\rho, \mu}(t)))),1\right) \nonumber\\
	&\stackrel{d}{=}  N\left(S_{\alpha'}(\lambda^{\alpha'}E_{\beta'}((M_{\alpha,\rho, \mu}(t)))),1\right) 
	\stackrel{d}{=}  N\left(\lambda S_{\alpha'}(E_{\beta'}((M_{\alpha,\rho, \mu}(t)))),1\right) 
	\stackrel{d}{=}  N\left( S_{\alpha'}(E_{\beta'}((M_{\alpha,\rho, \mu}(t)))),\lambda\right).
	\end{align}
	It may be noticed that (\ref{nhstfpp1}) with (\ref{c1}) may be viewed as the STFPP time changed by the ML L\'{e}vy subordinator.
\end{example}
\section{Simulation}
In this section, we reproduce some algorithms to simulate the sample paths of stable process, gamma process, ML L\'{e}vy process, and FPP. Using these, we present an algorithm to simulate the sample paths for the GFNBP.\\

\noindent {\bf Algorithm 1} (Simulation Algorithm for the Fractional Poisson Process).\\
Upto a fixed time T, the following algorithm gives the $n$ number of events $N_\beta(t)$ (see \cite{cahoy2010parameter, maheshwari2019fractional}).\\
(a) Fix the parameters $\lambda$ and $0< \beta < 1.$\\
(b) Set the initialization as $n=0$ and $t=0$.\\
(c) While $t < T$, generate three independent uniform random variables $U_i \sim U(0,1), i=1,2,3.$ Compute the increment as
	\begin{equation*}
	\Delta t = \frac{|\ln U_1|^{1/\beta}}{\lambda^{1/\beta}} \frac{\sin(\beta \pi U_2)(\sin(1-\beta)\pi U_2)^{1/{\beta -1}}}{(\sin{(\pi U_2)}^{1/\beta})|\ln U_3|^{1/{\beta-1}}}.
	\end{equation*}
(d) Update the increment as $t=t+\Delta t$ and $n=n+1.$\\
(e) Next t.\\

\noindent {\bf Algorithm 2} (Simulation Algorithm for the Stable Subordinator with $0< \alpha < 1)$.\\
(a) Choose $n$ time points $t_1, t_2,\dots, t_n.$ \\
(b) Generate an array of uniform random variables $U_i \sim U[0,\pi]$ and an array of exponential random variables $V_i \sim$ Exp$(1)$ for $i=1,2,\dots,n$.\\
(c) For 1 $\leq i \leq n,$ compute the increments as (see \cite{beghin2018space})
	\begin{equation*}
	\Delta S_\alpha^i t = S_\alpha(t_i) - S_\alpha(t_{i-1}) = (t_i - t_{i-1})^{1/\alpha} \frac{\sin(\alpha U_i)(\sin(1-\alpha) U_i)^{(1-\alpha)/\alpha}}{(\sin{( U_i)}^{1/\alpha})V_i^{(1-\alpha)/{\alpha}}}.
	\end{equation*}
(d) The sample path of $S_\alpha(t)$ at $t_i$ is $S_\alpha(t_i) = \sum_{j=1}^{i} 	\Delta S_\alpha^j.$\\

\noindent {\bf Algorithm 3} (Simulation Algorithm for the Gamma Subordinator).\\
(a) Fix the parameters $\rho$ and $\mu$.\\
(b) For a fixed time interval, choose equally spaced time points $t_i, i=1,2,\dots,n.$\\
(c) Generate $n$ independent gamma random variables $\mathcal{Q}_i \sim G(\mu, \rho t) $ for $i=1,2,\dots,n$ using the gamma sequential sampling technique (see \cite{avramidis2003efficient}).\\
(d) The sample path of gamma process $\Gamma(t)$ at $t_i$ is $\Gamma(ih) = \Gamma(t_i) = \sum_{j=1}^{i} \mathcal{Q}_j$ with $\mathcal{Q}_0 =0,$ where $h = t_2 - t-1$.\\

\noindent {\bf Algorithm 4} (Simulation Algorithm for the ML L\'{e}vy Subordinator).\\
(a) Fix the parameters $\rho$ and $\mu$.\\
(b) For a fixed time interval $[0,T]$, choose equally spaced time points $t_1 = t/n,\dots, (n-1)t/n=t_{n-1},T=t_n.$\\
(c) Generate a vector of size $n$ of the gamma random variables as $\mathcal{Q} = (\mathcal{Q}_1, \mathcal{Q}_2,\dots, \mathcal{Q}_n)$ such that $\mathcal{Q}_i \sim G(\mu, \rho(t_i - t_{i-1})).$\\
(d) Generate a vector of size $n$ of the $\alpha$-stable random variables as $S = (S_1, S_2,\dots, S_n)$ using Algorithm 2.\\
(e) Compute the increments of ML L\'{e}vy process via the self similar approach (see \cite{kumar2019tempered}) as $Y = (\mathcal{Q}_1^{1/{\mu}}S_1,\mathcal{Q}_2^{1/{\mu}}S_2,\dots,\mathcal{Q}_n^{1/{\mu}}S_n).$ \\
(f) Let $M_i =\sum_{ j=1}^{i}Y_j.$ Then, $M_1, M_2,\dots,M_n$ becomes the $n$ simulated values of the ML L\'{e}vy subordinator.\\

We next produce the algorithm to simulate the sample paths of the GFNBP using the above discussed algorithms.\\

\noindent {\bf Algorithm 5} (Simulation Algorithm for the GFNBP).\\
(a) Fix the parameters $\lambda$ and $\beta$ for the fractional Poisson process.\\
(b) For a fixed time interval $[0,T]$, choose equally spaced $(n+1)$ time points $t_0=0,t_1 = t/n,\dots, (n-1)t/n=t_{n-1},T=t_n$ with $h = t_2 -t_1$.\\
(c) Simulate the ML Levy subordinator $M(t_i)$ for $i=1,2,\dots,n$ using Algorithm 4.\\
(d) Using Algorithm 1, compute the number of events (arrivals) of the process $N_\beta(M(t_i))$ for $i = 1,2,\dots, n.$\\

\begin{figure}[h]
	\centering
	\subfloat{\includegraphics[width=9cm,height=6.25cm,keepaspectratio]{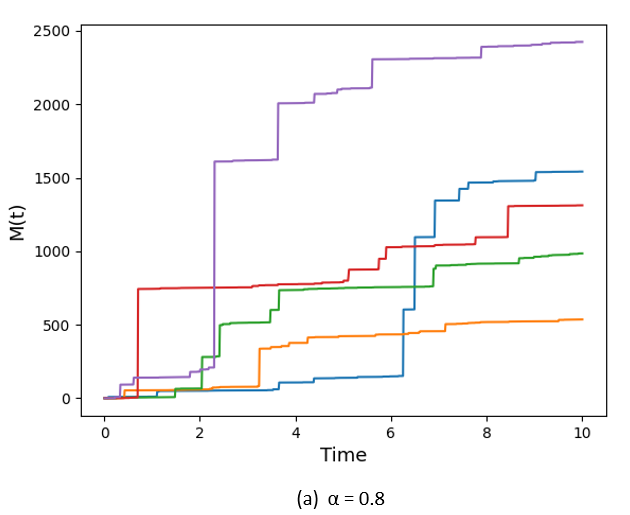}}
	\subfloat{\includegraphics[width=9cm,height=6.25cm,keepaspectratio]{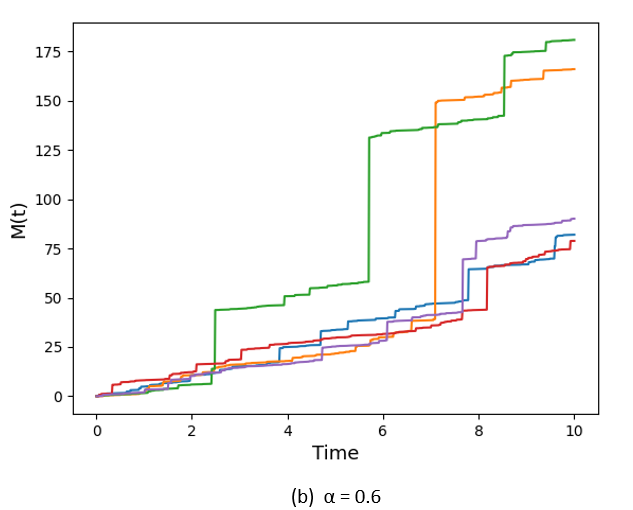}}
	\caption{The sample paths for the ML Levy Subordinator.}
	\label{Fig.1}
\end{figure}

\begin{figure}[ht]
	\centering
	\subfloat{\includegraphics[width=9cm,height=6.25cm,keepaspectratio]{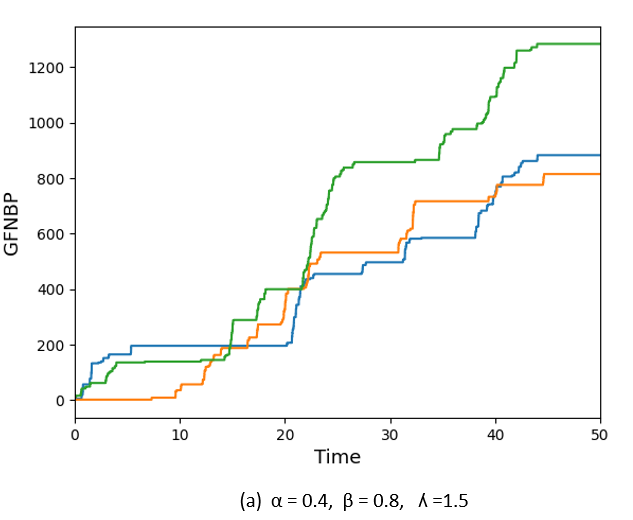}}
	\subfloat{\includegraphics[width=9cm,height=6.25cm,keepaspectratio]{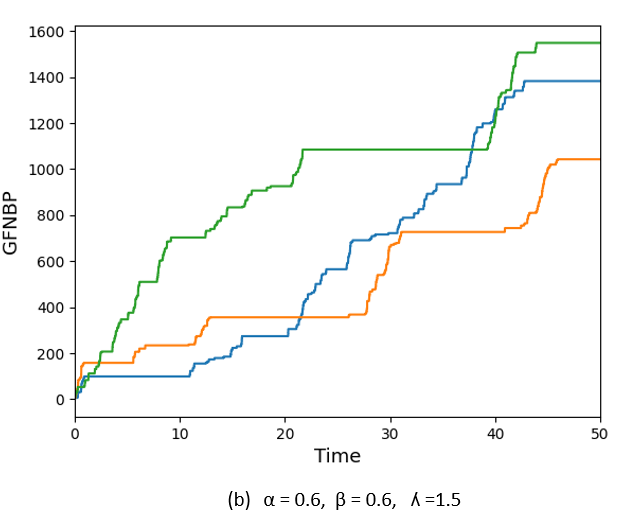}}
	\caption{The sample paths for the GFNBP.}
	\label{Fig.2}
\end{figure}

Based on the algorithm, the simulated paths for the ML L\'{e}vy subordinator and GFNBP are presented in Fig. \ref{Fig.1}  and Fig. \ref{Fig.2}, respectively.
\\
\section*{Declaration of competing interest}
The authors declare that they have no known competing financial interests.
\section*{Data availability}
No data was used for the research described in the article.
\bibliographystyle{apalike}
\bibliography{b4}

\begin{thebibliography}{33}
\providecommand{\natexlab}[1]{#1}
\providecommand{\url}[1]{\texttt{#1}}
\expandafter\ifx\csname urlstyle\endcsname\relax
  \providecommand{\doi}[1]{doi: #1}\else
  \providecommand{\doi}{doi: \begingroup \urlstyle{rm}\Url}\fi

\bibitem[Avramidis et~al.(2003)Avramidis, L~Ecuyer, Tremblay,
  et~al.]{avramidis2003efficient}
A.~N. Avramidis, P.~L~Ecuyer, P.-A. Tremblay, et~al.
\newblock Efficient simulation of gamma and variance-gamma processes.
\newblock In \emph{Winter Simulation Conference}, volume~1, pages 319--326,
  2003.

\bibitem[Barndorff-Nielsen(2000)]{barndorff2000probability}
O.~E. Barndorff-Nielsen.
\newblock \emph{Probability densities and Levy densities}.
\newblock University of Aarhus. Centre for Mathematical Physics and
  Stochastics~…, 2000.

\bibitem[Beghin(2015)]{beghin2015fractional}
L.~Beghin.
\newblock Fractional gamma and gamma-subordinated processes.
\newblock \emph{Stochastic Analysis and Applications}, 33\penalty0
  (5):\penalty0 903--926, 2015.

\bibitem[Beghin and Macci(2014)]{beghin2014fractional}
L.~Beghin and C.~Macci.
\newblock Fractional discrete processes: compound and mixed poisson
  representations.
\newblock \emph{Journal of Applied Probability}, 51\penalty0 (1):\penalty0
  19--36, 2014.

\bibitem[Beghin and Vellaisamy(2018)]{beghin2018space}
L.~Beghin and P.~Vellaisamy.
\newblock Space-fractional versions of the negative binomial and polya-type
  processes.
\newblock \emph{Methodology and Computing in Applied Probability}, 20\penalty0
  (2):\penalty0 463--485, 2018.

\bibitem[Biard and Saussereau(2014)]{biard2014fractional}
R.~Biard and B.~Saussereau.
\newblock Fractional poisson process: long-range dependence and applications in
  ruin theory.
\newblock \emph{Journal of Applied Probability}, 51\penalty0 (3):\penalty0
  727--740, 2014.

\bibitem[Cahoy et~al.(2010)Cahoy, Uchaikin, and Woyczynski]{cahoy2010parameter}
D.~O. Cahoy, V.~V. Uchaikin, and W.~A. Woyczynski.
\newblock Parameter estimation for fractional poisson processes.
\newblock \emph{Journal of Statistical Planning and Inference}, 140\penalty0
  (11):\penalty0 3106--3120, 2010.

\bibitem[Cook and Wei(2003)]{cook2003conditional}
R.~J. Cook and W.~Wei.
\newblock Conditional analysis of mixed poisson processes with baseline counts:
  implications for trial design and analysis.
\newblock \emph{Biostatistics}, 4\penalty0 (3):\penalty0 479--494, 2003.

\bibitem[Cox and Lewis(1966)]{cox1966statistical}
D.~R. Cox and P.~A. Lewis.
\newblock The statistical analysis of series of events.
\newblock 1966.

\bibitem[Gorenflo and Mainardi(2015)]{gorenflo2015fractional}
R.~Gorenflo and F.~Mainardi.
\newblock On the fractional poisson process and the discretized stable
  subordinator.
\newblock \emph{Axioms}, 4\penalty0 (3):\penalty0 321--344, 2015.

\bibitem[Grandell(1997)]{grandell1997mixed}
J.~Grandell.
\newblock \emph{Mixed poisson processes}, volume~77.
\newblock CRC Press, 1997.

\bibitem[Guler~Dincer et~al.(2022)Guler~Dincer, Demir, and
  Yal{\c{c}}in]{guler2022forecasting}
N.~Guler~Dincer, S.~Demir, and M.~O. Yal{\c{c}}in.
\newblock Forecasting covid19 reliability of the countries by using
  non-homogeneous poisson process models.
\newblock \emph{New Generation Computing}, pages 1--22, 2022.

\bibitem[Ken-Iti(1999)]{ken1999levy}
S.~Ken-Iti.
\newblock \emph{L{\'e}vy processes and infinitely divisible distributions}.
\newblock Cambridge university press, 1999.

\bibitem[Kilbas et~al.(2002)Kilbas, Saigo, and Trujillo]{kilbas2002generalized}
A.~A. Kilbas, M.~Saigo, and J.~J. Trujillo.
\newblock On the generalized wright function.
\newblock \emph{Fractional Calculus and Applied Analysis}, 5\penalty0
  (4):\penalty0 437--460, 2002.

\bibitem[Kumar et~al.(2019)Kumar, Upadhye, Wylomanska, and
  Gajda]{kumar2019tempered}
A.~Kumar, N.~Upadhye, A.~Wylomanska, and J.~Gajda.
\newblock Tempered mittag-leffler levy processes.
\newblock \emph{Communications in Statistics-Theory and Methods}, 48\penalty0
  (2):\penalty0 396--411, 2019.

\bibitem[Kumar et~al.(2020)Kumar, Leonenko, and Pichler]{kumar2020fractional}
A.~Kumar, N.~Leonenko, and A.~Pichler.
\newblock Fractional risk process in insurance.
\newblock \emph{Mathematics and Financial Economics}, 14:\penalty0 43--65,
  2020.

\bibitem[Laskin(2003)]{laskin2003fractional}
N.~Laskin.
\newblock Fractional poisson process.
\newblock \emph{Communications in Nonlinear Science and Numerical Simulation},
  8\penalty0 (3-4):\penalty0 201--213, 2003.

\bibitem[Laskin(2009)]{laskin2009some}
N.~Laskin.
\newblock Some applications of the fractional poisson probability distribution.
\newblock \emph{Journal of Mathematical Physics}, 50\penalty0 (11):\penalty0
  113513, 2009.

\bibitem[Maheshwari(2023)]{maheshwari2023tempered}
A.~Maheshwari.
\newblock Tempered space fractional negative binomial process.
\newblock \emph{Statistics \& Probability Letters}, 196:\penalty0 109799, 2023.

\bibitem[Maheshwari and Vellaisamy(2016)]{maheshwari2016long}
A.~Maheshwari and P.~Vellaisamy.
\newblock On the long-range dependence of fractional poisson and negative
  binomial processes.
\newblock \emph{Journal of Applied Probability}, 53\penalty0 (4):\penalty0
  989--1000, 2016.

\bibitem[Maheshwari and
  Vellaisamy(2019{\natexlab{a}})]{maheshwari2019fractional}
A.~Maheshwari and P.~Vellaisamy.
\newblock Fractional poisson process time-changed by l{\'e}vy subordinator and
  its inverse.
\newblock \emph{Journal of Theoretical Probability}, 32\penalty0 (3):\penalty0
  1278--1305, 2019{\natexlab{a}}.

\bibitem[Maheshwari and Vellaisamy(2019{\natexlab{b}})]{maheshwari2019non}
A.~Maheshwari and P.~Vellaisamy.
\newblock Non-homogeneous space-time fractional poisson processes.
\newblock \emph{Stochastic Analysis and Applications}, 37\penalty0
  (2):\penalty0 137--154, 2019{\natexlab{b}}.

\bibitem[Mathai et~al.(2009)Mathai, Saxena, and Haubold]{mathai2009h}
A.~M. Mathai, R.~K. Saxena, and H.~J. Haubold.
\newblock \emph{The H-function: theory and applications}.
\newblock Springer Science \& Business Media, 2009.

\bibitem[Meerschaert et~al.(2011)Meerschaert, Nane, and
  Vellaisamy]{meerschaert2011fractional}
M.~Meerschaert, E.~Nane, and P.~Vellaisamy.
\newblock The fractional poisson process and the inverse stable subordinator.
\newblock \emph{Electronic Journal of Probability}, 16:\penalty0 1600--1620,
  2011.

\bibitem[Meerschaert and Scheffler(2004)]{meerschaert2004limit}
M.~M. Meerschaert and H.-P. Scheffler.
\newblock Limit theorems for continuous-time random walks with infinite mean
  waiting times.
\newblock \emph{Journal of applied probability}, 41\penalty0 (3):\penalty0
  623--638, 2004.

\bibitem[Meerschaert and Straka(2013)]{meerschaert2013inverse}
M.~M. Meerschaert and P.~Straka.
\newblock Inverse stable subordinators.
\newblock \emph{Mathematical modelling of natural phenomena}, 8\penalty0
  (2):\penalty0 1--16, 2013.

\bibitem[Orsingher and Polito(2012)]{orsingher2012space}
E.~Orsingher and F.~Polito.
\newblock The space-fractional poisson process.
\newblock \emph{Statistics \& Probability Letters}, 82\penalty0 (4):\penalty0
  852--858, 2012.

\bibitem[Podlubny(1999)]{podlubny1999introduction}
I.~Podlubny.
\newblock An introduction to fractional derivatives, fractional differential
  equations, to methods of their solution and some of their applications.
\newblock \emph{Math. Sci. Eng}, 198:\penalty0 340, 1999.

\bibitem[Polito and Scalas(2016)]{polito2016generalization}
F.~Polito and E.~Scalas.
\newblock A generalization of the space-fractional poisson process and its
  connection to some l{\'e}vy processes.
\newblock 2016.

\bibitem[Prabhakar et~al.(1971)]{prabhakar1971singular}
T.~R. Prabhakar et~al.
\newblock A singular integral equation with a generalized mittag-leffler
  function in the kernel.
\newblock \emph{Yokohama math. J}, 19\penalty0 (1):\penalty0 7--15, 1971.

\bibitem[Shukla and Prajapati(2007)]{shukla2007generalization}
A.~Shukla and J.~Prajapati.
\newblock On a generalization of mittag-leffler function and its properties.
\newblock \emph{Journal of mathematical analysis and applications},
  336\penalty0 (2):\penalty0 797--811, 2007.

\bibitem[Timmermann and Nowak(1999)]{timmermann1999multiscale}
K.~E. Timmermann and R.~D. Nowak.
\newblock Multiscale modeling and estimation of poisson processes with
  application to photon-limited imaging.
\newblock \emph{IEEE Transactions on Information Theory}, 45\penalty0
  (3):\penalty0 846--862, 1999.

\bibitem[Vellaisamy and Maheshwari(2018)]{samy2018fractional}
P.~Vellaisamy and A.~Maheshwari.
\newblock Fractional negative binomial and polya processes.
\newblock \emph{Probability and Mathematical Statistics}, 38\penalty0
  (1):\penalty0 77--101, 2018.

\end{thebibliography}











\end{document}